# POWER LAWS FOR FAMILY SIZES IN A DUPLICATION MODEL

By Rick Durrett[1] and Jason Schweinsberg[2]

*Cornell University and University of California, San Diego*

Qian, Luscombe and Gerstein [*J. Molecular Biol.* **313** (2001) 673–681] introduced a model of the diversification of protein folds in a genome that we may formulate as follows. Consider a multitype Yule process starting with one individual in which there are no deaths and each individual gives birth to a new individual at rate 1. When a new individual is born, it has the same type as its parent with probability $1 - r$ and is a new type, different from all previously observed types, with probability $r$. We refer to individuals with the same type as families and provide an approximation to the joint distribution of family sizes when the population size reaches $N$. We also show that if $1 \ll S \ll N^{1-r}$, then the number of families of size at least $S$ is approximately $CNS^{-1/(1-r)}$, while if $N^{1-r} \ll S$ the distribution decays more rapidly than any power.

**1. Introduction.** Genome sequencing of various species has shown that gene and protein-fold family sizes have a power-law distribution. Huynen and van Nimwegen [19] studied six bacteria, two Archea and yeast. Li, Gu, Wang and Nekrutenko [28] and later Gu, Cavalcanti, Chen, Bouman and Li [17] analyzed the genomes of yeast, the nematode *C. elegans*, fruit fly (*Drosophila melanogaster*) and human. There have been several models advanced to explain this phenomenon. Rzhetsky and Gomez [33] and Karev, Wolf, Rzhetsky, Berezov and Koonin [23] (see also [25]) introduced a birth and death model in which, when there are $i$ individuals in a family, a birth occurs at rate $\lambda_i$ and a death occurs at rate $\delta_i$. They proved, as most readers of this journal can easily verify, that if the birth rates are second-ordered balanced, that is,

$$\lambda_{i-1}/\delta_i = 1 - a/i + O(1/i^2)$$

Received May 2004; revised November 2004.

[1]Supported in part by NSF Grants 0202935 from the probability program and 0201037 from a joint DMS/NIGMS initiative to support research in mathematical biology.

[2]Supported by an NSF Postdoctoral Fellowship.

for some $a > 0$, then the stationary distribution of the family size is asymptotically $Ci^{-a}$. See the Appendix of [23] or Example 3.6 on page 297 of [13] for more details.

Qian, Luscombe and Gerstein [32] introduced an alternative model that we will study in detail here. Consider a continuous-time Yule process with infinitely many types. At time zero, a single individual of type 1 is born. No individuals die, and each individual independently gives birth to a new individual at rate 1. When a new individual is born, it has the same type as its parent with probability $1 - r$, where $0 < r < 1$. With probability $r$, the new individual has a type which is different from all previously observed types. If the $k$th individual born has a different type from its parent, we say that it has type $k$. Note that, as a consequence of this choice of labeling, there are always type-1 individuals, but for $k \geq 2$, with probability $1 - r$ there are never any individuals of type $k$.

In this model one can think of the new types as resulting from mutations, where $r$ is the probability of mutation. Alternatively, one could think of a Yule process with immigration in which each individual gives birth at rate $1 - r$ and new immigrants arrive at rate $r$ times the current population size. We refer to individuals with the same type as families. The goal of this paper is to study the distribution of the family sizes at the time when the population size reaches $N$.

1.1. *Approximation to the family-size distribution.* Let $T_N$ be the time that the population size reaches $N$. Let $R_{k,N}$ be the number of individuals of type $k$ at time $T_N$. Let $X_{k,N}$ be the fraction of individuals at time $T_N$ whose type is in $\{1, \ldots, k\}$. Let $V_{k,N}$ be the fraction of individuals at time $T_N$, among those whose type is in $\{1, \ldots, k\}$, that are of type $k$. This means that the fraction of individuals at time $T_N$ that are of type $k$ is $V_{k,N} X_{k,N}$ and the number of individuals of type $k$ at time $T_N$ is $R_{k,N} = N V_{k,N} X_{k,N}$. Note that $X_{N,N} = 1$ and for $k = 1, \ldots, N - 1$, we have

$$(1.1) \qquad X_{k,N} = \prod_{j=k+1}^{N} (1 - V_{j,N}).$$

The following proposition follows from well-known connections between Yule processes and Pólya urns. We review these connections and prove this proposition in Section 2.

PROPOSITION 1.1. *For each positive integer $k$, the limit*

$$W_k = \lim_{N \to \infty} V_{k,N}$$

*exists a.s. The random variables $W_1, W_2, \ldots$ are independent. We have $W_1 = 1$ a.s. Furthermore, $P(W_k > 0) = r$ for all $k \geq 2$ and conditional on the event that $W_k > 0$, the distribution of $W_k$ is* Beta$(1, k-1)$.

Let $Y_{N,N} = 1$ and, for $k = 1, \ldots, N - 1$, let

$$(1.2) \qquad Y_{k,N} = \prod_{j=k+1}^{N} (1 - W_j).$$

Let $\Delta = \{(x_i)_{i=1}^{\infty} : 0 \leq x_i \leq 1 \text{ for all } i \text{ and } \sum_{i=1}^{\infty} x_i = 1\}$. Note that the sequence

$$(N^{-1} R_{k,N})_{k=1}^{\infty} = (V_{k,N} X_{k,N})_{k=1}^{\infty},$$

whose $k$th term is the fraction of the population having type $k$ at time $T_N$, is in $\Delta$. Proposition 1.1 and equations (1.1) and (1.2) suggest that, for large $N$, the distribution of this sequence can be approximated by $Q_{r,N}$, which we define to be the distribution of the sequence in $\Delta$ whose first term is $Y_{1,N}$, whose $k$th term is $W_k Y_{k,N}$ for $2 \leq k \leq N$ and whose $k$th term is zero for $k > N$. Theorem 1.2 below uses the coupling of the $X_{k,N}$ and $Y_{k,N}$ given above to show that the distribution of $(N^{-1} R_{k,N})_{k=1}^{\infty}$ can be approximated by $Q_{r,N}$ to within an error of $O(N^{-1/2})$. We prove this result in Section 3.

THEOREM 1.2. *We have* $E[\max_{1 \leq k \leq N} |X_{k,N} - Y_{k,N}|] \leq \frac{5}{\sqrt{N}}$.

The distributions $Q_{r,N}$ first arose in the work of Durrett and Schweinsberg [15] and Schweinsberg and Durrett [34], who studied the effect of beneficial mutations on the genealogy of a population. The distributions $Q_{r,N}$ arose in that context because, shortly after a beneficial mutation, the number of individuals with the beneficial gene behaves like a supercritical branching process, which means that the number with descendants surviving a long time into the future behaves like a Yule process. In this setting, $r$ is the rate of recombination, and individuals descended from a lineage with a recombination get traced back to a different ancestor than other individuals, just as individuals descended from an individual with a mutation in the present model are of a different type than the others. Schweinsberg and Durrett's [34] approximation had an error of $O((\log N)^{-2})$ because of deaths and other complexities in the model, but Theorem 1.2 shows that the distributions $Q_{r,N}$ give a much more accurate approximation to the family-size distribution in the simpler model studied here. We note also that here it is assumed that $r$ is fixed, whereas Schweinsberg and Durrett [34] considered $r$ to be $O(1/(\log N))$.

1.2. *A power law for the number of families of moderate size.* Let $F_{S,N}$ denote the number of families at time $T_N$ whose size is at least $S$. Define

$$(1.3) \qquad g(S) = r \Gamma \left( \frac{2 - r}{1 - r} \right) N S^{-1/(1-r)}.$$

The theorem below, which is proved in Section 4, shows that if $1 \ll S \ll N^{1-r}$, then $g(S)$ provides a good approximation to the number of families of size at least $S$, in the sense that $|F_{S,N} - g(S)|/g(S) \to 0$ as $N \to \infty$.

THEOREM 1.3.   *There are constants $0 < C_1, C_2 < \infty$ so that*

$$E[|F_{S,N} - g(S)|] \leq C_1 g(S)[S^{-1/5} + (NS^{-1/(1-r)})^{-1/5}] + C_2.$$

Note that $S^{-1/5}$ and $(NS^{-1/(1-r)})^{-1/5}$ are both small and $g(S)$ is large when $1 \ll S \ll N^{1-r}$.

Theorem 1.3 confirms Qian, Luscombe and Gerstein's [32] power law but it also conflicts with their results. Since they considered the number of folds that occur exactly $V$ times rather than at least $V$ times, it follows from differentiating the right-hand side of (1.3) that for large $N$ we would expect a decay with the power $b = 1 + 1/(1 - r)$. This quantity is always larger than 2, while they observed powers $b$ between 0.9 and 1.2 for eukaryotes and between 1.2 and 1.8 for prokaryotes. Despite this discrepancy, they were able to fit their model by starting the process at time zero with $N_0 > 1$ families. For example, for *Haemophilus influenzae* they took $r = 0.3$, $N_0 = 90$, and ran the process for 1249 generations. For *C. elegans* they took $r = 0.018$, $N_0 = 280$ and ran for 18,482 generations.

Figure 1 shows one simulation of the system with $r = 0.018$, $N_0 = 1$, and $N = 20,000$. In contrast to biologists who do a log-log plot of the number of gene families of size $k$ (see, e.g., Figure 1 in [18], or Figure 8 in [23]), we look at the tail of the distribution and plot the log of the family size on the $x$-axis and the log of the number of families of at least that size on the $y$-axis. The curve fit by Karev et al. [23] has asymptotic power 1.9 in contrast to the 2.018 that comes from our formula, but note that the straight line fitted to our simulation of the distribution function has slope 0.91. Figure 2 shows the average of 10,000 simulations of the process with the *C. elegans* parameters. The straight line shows that Theorem 1.3 very accurately predicts the expected number of families until the log of the family size is 4. This simulation also shows that the power law breaks down when $S \gg N^{1-r}$, which motivates our next topic.

1.3. *Sizes of the largest families.*   Recall that $R_{k,N}$ is the number of individuals of type $k$ at time $T_N$. Proposition 1.4 below identifies the limiting distribution of the size of the large families.

PROPOSITION 1.4.   *For each positive integer $k$, the limit*

$$Z_k = \lim_{N \to \infty} N^{r-1} R_{k,N}$$

*exists almost surely. The distribution of $Z_1$ is the Mittag–Leffler distribution with parameter $1 - r$, which has density*

$$(1.4) \quad g(x) = \frac{1}{\pi(1-r)} \sum_{k=0}^{\infty} \frac{(-1)^{k+1}}{k!} \sin(\pi\alpha k)\Gamma((1-r)k+1)x^{k-1}, \qquad x > 0.$$

*For $k \geq 2$, conditional on the event that the $k$th individual born has type $k$, the distribution of $Z_k$ is the same as the distribution of $MB_k^{1-r}$, where $B_k$ has the Beta$(1, k-1)$ distribution, $M$ has the Mittag–Leffler distribution with parameter $1 - r$ and $M$ and $B_k$ are independent.*

The fact that $Z_1$ has the Mittag–Leffler distribution was first proved by Angerer [2], who was motivated by the study of bacterial populations. He considered a model that is equivalent to our model, except that he referred to our type-1 individuals as nonmutant cells, and individuals of all other types as mutant cells. Theorem 6.1 in [2] gives the Mittag–Leffler limit when the probability of mutation is a fixed constant. See also Theorems 1.7 and 1.8 of [21], where the Mittag–Leffler distribution arises as a limiting distribution in an urn model that is closely related to our model. We mention another proof of the Mittag–Leffler limit at the end of Section 1.4, and we prove Proposition 1.4 for $k \geq 2$ in Section 5.

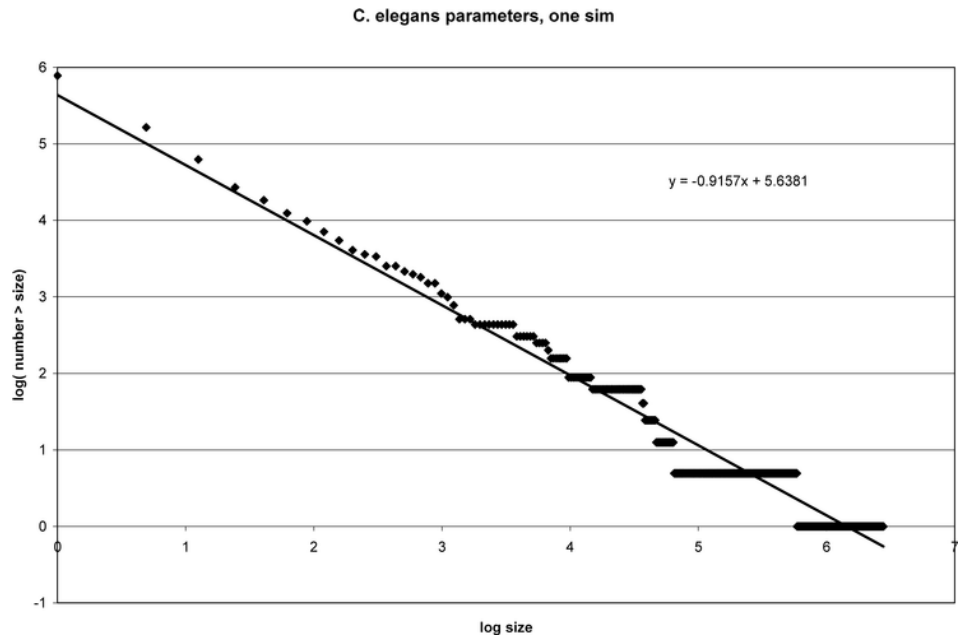

FIG. 1. *One simulation of the duplication model with C. elegans parameters. $r = 0.018$ and $N = 20,000$.*

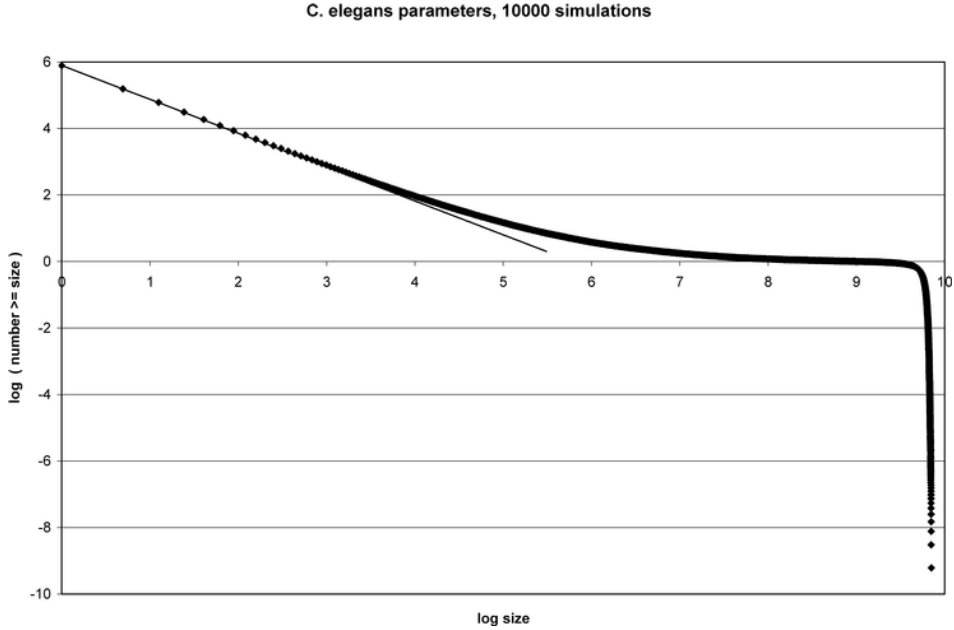

Fig. 2. *Average of* 10,000 *simulations of the duplication model with C. elegans parameters.* $r = 0.018$ *and* $N = 20,000$. *Straight line is the prediction of Theorem* 1.3.

The moments of $M$ are given by $E[M^m] = \Gamma(m+1)/\Gamma(m(1-r)+1)$ for $m > 0$ (see Section 0.5 of [30]). Also, we have $E[B_k^m] = \Gamma(m+1)\Gamma(k)/\Gamma(m+k)$ for $m > 0$. Since $P(Z_k > 0) = r$ when $k \geq 2$, it follows that for $k \geq 2$ and $m > 0$ we have

$$E[Z_k^m] = \frac{\Gamma(m+1)}{\Gamma(m(1-r)+1)} \cdot \frac{r\Gamma(m(1-r)+1)\Gamma(k)}{\Gamma(m(1-r)+k)} = \frac{r\Gamma(m+1)\Gamma(k)}{\Gamma(m(1-r)+k)}.$$

The next result, which is proved in Section 5, proves what was observed in the simulation. The expected number of families of size at least $xN^{1-r}$ decays faster than any power of $x$. Indeed, it decays faster than exponentially in $x$, and the decay is fastest when $r$ is small.

PROPOSITION 1.5. *There exist constants $C_1$ and $C_2$ such that for all $x \geq 1$, we have*

$$\lim_{N \to \infty} \sum_{k=1}^{N} P(R_{k,N} > xN^{1-r}) \leq C_1 e^{-C_2 x^{1/r}}.$$

1.4. *A new Chinese restaurant.* Our model has a close relation to a construction called the "Chinese restaurant process," which was first proposed by Dubins and Pitman. We describe here a two-parameter version of the process, which is discussed in Pitman [29, 30]. Suppose $0 \leq \alpha < 1$ and $\theta > -\alpha$. Consider a restaurant with infinitely many tables, each with an unbounded number of seats. The first customer sits at table 1. Suppose, for some $n \geq 1$, that after $n$ customers have been seated, there are $k$ occupied tables, with $n_i$ customers at the $i$th table, so that $n_1 + \cdots + n_k = n$. Then, the $(n+1)$st customer sits at table $i$ with probability $(n_i - \alpha)/(n + \theta)$ and sits at an unoccupied table, which we call the $(k+1)$st table, with probability $(\theta + k\alpha)/(n + \theta)$.

For any $N$, the Chinese restaurant process gives rise to a random partition $\Pi_N$ of $\{1, \ldots, N\}$, where $i$ and $j$ are in the same block of $\Pi_N$ if and only if the $i$th and $j$th customers are seated at the same table. That is, the partition $\Pi_N$ consists of blocks $B_{1,N}, \ldots, B_{k,N}$, where $B_{j,N}$ consists of all integers $i$ between 1 and $N$ such that the $i$th customer is seated at the $j$th table. Let $|B_{j,N}|$ denote the number of the first $N$ customers at the $j$th table. Then (see [30]), the distribution of the $\Delta$-valued sequence $(N^{-1}|B_{1,N}|, N^{-1}|B_{2,N}|, \ldots)$ converges as $N \to \infty$ to the Poisson–Dirichlet distribution with parameters $(\alpha, \theta)$. This distribution is defined as follows. Let $(D_j)_{j=1}^{\infty}$ be a sequence of independent random variables such that $D_j$ has the Beta$(1 - \alpha, \theta + j\alpha)$ distribution. Then the sequence whose $k$th term is $D_k \prod_{j=1}^{k-1}(1 - D_j)$ has the Poisson–Dirichlet distribution with parameters $(\alpha, \theta)$. The Poisson–Dirichlet distributions were studied extensively by Pitman and Yor [31]. See also [30] and [4] for further applications of these distributions.

An important special case of the Chinese restaurant process is when $\alpha = 0$. Then, we may assume that the $(n+1)$st customer sits at a new table with probability $\theta/(n + \theta)$ and otherwise chooses one of the previous $n$ customers at random and sits at that person's table. In this case, if $\pi$ is a partition of $\{1, \ldots, N\}$ with $k$ blocks of sizes $n_1, \ldots, n_k$, one can check that

$$(1.5) \qquad P(\Pi_N = \pi) = \frac{\theta^{k-1}}{(1 + \theta)(2 + \theta) \cdots (N - 1 + \theta)} \prod_{i=1}^{k} (n_i - 1)!.$$

This leads to the famous Ewens sampling formula [16]. The Ewens sampling formula describes the family-size distribution in a Yule process with immigration when immigration occurs at constant rate $\theta$. When there are $n$ individuals in the Yule process, they are each splitting at rate 1 and immigration occurs at rate $\theta$, so the probability that the $(n + 1)$st individual starts a new family is $\theta/(n + \theta)$. For another application of the Ewens sampling formula, consider a population in which each lineage experiences mutation at rate $\theta/2$ and whose ancestral structure is given by Kingman's coalescent (see [24]), meaning that each pair of lineages merges at rate 1. Working backward in time, when there are $n + 1$ lineages, coalescence occurs at rate $n(n+1)/2$ while mutations occur at rate $\theta(n+1)/2$. Consequently, the probability of having mutation before coalescence is $\theta/(n + \theta)$. Because Kingman's coalescent is a good approximation to the genealogy in populations of fixed size, the Ewens sampling formula is a standard model for gene frequencies in populations of fixed size. However, this model does not lead to the power-law behavior that has been observed in some data.

Note that our model can be viewed as a variation of the Chinese restaurant process in which the $(n+1)$st customer sits at a new table with constant probability $r$, rather than with probability $\theta/(n+\theta)$, and otherwise picks one of the previous $n$ customers at random and sits at that person's table. One can define the random partition $\Theta_N$ of $\{1, \ldots, N\}$ such that $i$ and $j$ are in the same block if and only if the $i$th and $j$th customers are seated at the same table. In our branching process interpretation, this means that the $i$th and $j$th individuals born have the same type, so the family sizes in our model correspond to block sizes of $\Theta_N$. It is straightforward to derive an analog of the Ewens sampling formula in this case. If $\pi$ is a partition of $\{1, \ldots, N\}$ into $k$ blocks of sizes $n_1, \ldots, n_k$, and if $a_1 < a_2 < \cdots < a_k$ are the first integers in these blocks, then

$$P(\Theta_N = \pi) = \frac{r^{k-1}(1-r)^{N-k}}{(N-1)!} \left[ \prod_{i=1}^{k} (n_i - 1)! \right] \prod_{j=2}^{k} (a_j - 1).$$

This formula depends on $a_2, \ldots, a_k$ as well as the block sizes $n_1, \ldots, n_k$, so the random partition $\Theta_N$ is not exchangeable. Nevertheless, one can still look for approximations to the distribution of the block sizes. We see from Theorem 1.2 that the distributions $Q_{r,N}$ play the role of the Poisson–Dirichlet distributions in this model. Because the population size in Yule processes grows exponentially, these distributions provide a plausible model of gene frequencies in growing populations, and they do lead to power-law behavior, as shown in Theorem 1.3. Furthermore, the approximation error in Theorem 1.2 of $O(N^{-1/2})$ is the same order of magnitude as the error when the distributions of the block sizes of the partitions $\Pi_N$ above are approximated by the Poisson–Dirichlet distributions.

Finally, the Chinese restaurant process when $\alpha = 1 - r$ and $\theta = 0$ can be used to give another proof of Proposition 1.4 when $k = 1$. This argument was pointed out to us by Wolfgang Angerer, Anton Wakolbinger and a referee, and also appears implicitly in earlier unpublished notes of Jim Pitman. Given our multitype Yule process, we can obtain a Chinese restaurant process with $\alpha = 1 - r$ and $\theta = 0$ by saying that each individual born in the Yule process sits at the same table as its parent, unless it has type 1 in which case it starts a new table. Thus, the number of type-1 individuals is the number of occupied tables, so the Mittag–Leffler limit follows from Theorem 31 of [30]. See also Angerer and Wakolbinger [3].

1.5. *Connections with preferential attachment.* In our model, gene families grow at a rate proportional to their size. This is similar to the behavior of Barbási and Albert's [8] preferential attachment model in which one grows a graph by adding a vertex at each time and connecting that vertex to $m$ existing vertices chosen with probabilities proportional to their degrees. Through simulations and heuristic arguments, Barbási and Albert concluded that the fraction of vertices of degree $k$ converged to a limit $p_k \sim Ck^{-3}$. This result was later proved rigorously by Bollobás, Riordan, Spencer and Tusnády [11].

Fueled by the observation of power laws for degree distributions in the Internet, collaboration networks and even sexual relations in Sweden, this work touched off a flurry of activity. To remedy the difficulty that the power was always 3 in the Barbási–Albert model, Krapivsky, Redner and Leyvraz [26] introduced a model in which attachment to vertices of degree $i$ was proportional to $a + bi$, and were able to achieve any power in $(2, \infty)$. These results, published in *Physical Review Letters*, omit a few details, but work by Kumar, Raghavan, Rajagopalan, Sivakumar, Tompkins and Upfal [27] and Cooper and Freize [12] further generalizes the model and provides rigorous proofs of the power laws.

The preferential attachment models are different from ours because adding an edge changes the degree of two vertices. However, if one considers directed graphs and analyzes only the out-degree, then taking $\alpha = 1 - r$ and $\delta = 0$ in the Cooper–Freize model gives a model identical to ours and a power law that is proved in their Section 6.1. Later work of Bollobás, Borgs, Chayes and Riordan [10] investigates a directed graph model which contains our result as a special case and for which they derive a power law. We point out also that a construction similar to that given in Section 1.1 was developed by Berger, Borgs, Chayes and Saberi [9] in the context of preferential attachment graphs.

In addition to recent work, Simon [35] considered the following model of word usage in books, which he also applied to scientific publications, city sizes and income distribution. Let $X_i(t)$ be the number of words that have appeared exactly $i$ times in the first $t$ words. He assumed that (a) the probability that the $(t + 1)$st word is a word that has already appeared $i$ times is proportional to $iX_i(t)$; (b) there is a constant probability $\alpha$ that the $(t + 1)$st word is a word that has not appeared in the first $t$ words. This of course is exactly our model, but even this is not the earliest reference. It appeared in work of Yule [37] who considered a model of the number of species in a given genus. Both Yule [37] and Simon [35] argued that the model gives rise to power-law behavior. See Aldous [1] for a more recent account and a simple explanation for the power law.

While our model has been considered a number of times, our results are more precise. In most cases investigators have considered the limit of the fraction of vertices of degree $k$ for fixed $k$. Exceptions are Bollobás, Riordan,

Spencer and Tusnády [11] who were able to prove results for $k \leq N^{1/15}$ and Cooper and Freize [12] who could handle $k \leq N^{1/21}$. In contrast, our results hold for the entire range over which the power law is valid and show how the power law breaks down for larger values.

## 2. Branching processes and Pólya urns.

In this section we review some well-known connections between Pólya urns and continuous-time branching processes, which will be useful later in the paper. Athreya and Karlin [6] showed how to embed the urn process in a continuous-time branching process. This technique was reviewed in [7]. See [20] for a thorough survey of recent developments and generalizations.

Recall the following version of Pólya's urn model. Suppose we start with $a$ white balls and $b$ black balls in the urn. We then draw a ball at random from the urn. If the ball we draw is white, we return it to the urn and add an additional white ball to the urn. If the ball we draw is black, we return it to the urn and add another black ball. This process can be repeated indefinitely. To see the connection with branching processes, consider a two-type branching process in which there are no deaths and each individual gives birth at rate 1. If at some time there are $a$ individuals of type 1 and $b$ individuals of type 2, then the probability that the next individual born will have type 1 is $a/(a+b)$, which is the same as the probability that the next ball added to an urn with $a$ white balls and $b$ black balls will be white. It follows that the distribution of the number of type-1 individuals when the population size reaches $N$ is the same as the distribution of the number of white balls in the urn when the number of balls in the urn is $N$.

Let $\zeta_i = 1$ if the $i$th ball added to the urn is white, and let $\zeta_i = 0$ if the $i$th ball added to the urn is black. Fix a positive integer $N$. Let $S \subseteq \{1, \ldots, N\}$, and let $S^c = \{1, \ldots, N\} \setminus S$. Let $|S|$ denote the cardinality of $S$. It is easy to check that for $a, b \geq 1$,

$$
\begin{aligned}
(2.1) \quad & P(\zeta_i = 1 \text{ for } i \in S \text{ and } \zeta_i = 0 \text{ for } i \in S^c) \\
& = \frac{(a + |S| - 1)!(b + N - |S| - 1)!(a + b - 1)!}{(a-1)!(b-1)!(a+b+N-1)!}.
\end{aligned}
$$

Since the right-hand side of (2.1) depends only on $|S|$ and not on the particular elements of $S$, the sequence $(\zeta_i)_{i=1}^{\infty}$ is exchangeable. By de Finetti's theorem, there exists a probability measure $\mu$ on $[0, 1]$ such that for all $N$ and all $S \subseteq \{1, \ldots, N\}$, we have

$$
(2.2) \quad P(\zeta_i = 1 \text{ for } i \in S \text{ and } \zeta_i = 0 \text{ for } i \in S^c) = \int_0^1 x^{|S|}(1-x)^{N-|S|}\mu(dx),
$$

where $\mu$ is the distribution of $\lim_{N \to \infty} N^{-1}|\{i \leq N : \zeta_i = 1\}|$, the limiting fraction of white balls in the urn when we start with $a$ white balls and $b$ black balls. It follows from Theorem 1 in Section 9.1 of Chapter V of [7] that $\mu$ is the Beta$(a, b)$ distribution. One can also see this by checking that the right-hand sides of (2.1) and (2.2) agree in this case.

PROOF OF PROPOSITION 1.1. Clearly $W_1 = 1$ a.s. because $V_{1,N} = 1$ a.s. Assume now that $k \geq 2$. Let $S_k$ be the set of all $i$ such that the type of the $i$th individual born is in $\{1, \ldots, k\}$. Let $\mathcal{H}_k$ be the $\sigma$-field generated by the sets $S_k, S_{k+1}, \ldots$. Note that if $j > k$, then $V_{j,N}$ is $\mathcal{H}_k$-measurable for all $N$. Therefore, to prove the proposition, it suffices to show that, for all $k \geq 2$, the limit $W_k$ exists a.s. and satisfies the following conditions:

1. $P(W_k > 0) = r$.
2. The conditional distribution of $W_k$ given $W_k > 0$ is Beta$(1, k-1)$.
3. $W_k$ is independent of $\mathcal{H}_k$.

Note that the third condition implies that $W_k$ is independent of $(W_j)_{j=k+1}^\infty$.

Enumerate the elements of $S_k$ as $i_1 < i_2 < i_3 < \cdots$. Define a sequence $(\zeta_j^{(k)})_{j=1}^\infty$ such that $\zeta_j^{(k)} = 1$ if the $i_j$th individual has type $k$ and $\zeta_i^{(k)} = 0$ otherwise. Note that $i_j = j$ for $j \leq k$. Also, $\zeta_j^{(k)} = 0$ for $j = 1, \ldots, k-1$. Recall from our conventions for labeling the types that if the $k$th individual to enter the population has a new type, then it has type $k$. Therefore, $\zeta_k^{(k)} = 1$ if and only if the $k$th individual has a new type, and whether or not this individual has a new type does not affect the births of individuals of types greater than $k$. Thus, $P(\zeta_k^{(k)} = 1|\mathcal{H}_k) = r$. If $\zeta_k^{(k)} = 0$, then clearly $W_k = 0$. Because of the connection between branching processes and Pólya urns, if $\zeta_k^{(k)} = 1$, then the sequence $(\zeta_{j+k}^{(k)})_{j=1}^\infty$ has the same distribution as the Pólya urn sequence $(\zeta_i)_{i=1}^\infty$ defined above when $a = 1$ and $b = k-1$. Furthermore, the values of $\zeta_j^{(k)}$ do not affect the births of individuals of types greater than $k$, so this relationship holds even after conditioning on $\mathcal{H}_k$. It follows that, conditional on $\zeta_k^{(k)} = 1$, the random variable $W_k$ has a Beta$(1, k-1)$ distribution and $W_k$ is independent of $\mathcal{H}_k$. □

Now, fix $N$ and to simplify notation, write $X_k$, $Y_k$, $V_k$ and $R_k$ for $X_{k,N}$, $Y_{k,N}$, $V_{k,N}$ and $R_{k,N}$, respectively. We will use this notation throughout the rest of the paper when the value of $N$ is clear from the context. Let $\mathcal{F}_k$ be the $\sigma$-field generated by the random variables $V_j$ and $W_j$ for $j \geq k+1$. It follows from (1.1) and (1.2) that $X_k$ and $Y_k$ are $\mathcal{F}_k$-measurable. Let $\mathcal{G}_k$ be the $\sigma$-field generated by the random variables $V_j$ for $j \geq k+1$ and $W_j$ for $j \geq k$.

We can write $W_k = \xi_k \tilde{W}_k$, where $\xi_k$ has a Bernoulli$(r)$ distribution and is independent of $\mathcal{F}_k$, and $\tilde{W}_k$ has a Beta$(1, k-1)$ distribution and is independent of $\xi_k$ and $\mathcal{F}_k$. Since $E[\tilde{W}_k] = 1/k$ and $E[\tilde{W}_k^2] = 2/[k(k+1)]$, we have $E[W_k|\mathcal{F}_k] = r/k$ and $E[W_k^2|\mathcal{F}_k] = 2r/[k(k+1)]$. Note that $V_k = 0$ whenever $W_k = 0$. On $\{W_k > 0\}$, define $\tilde{V}_k = R_k - 1$. Define $\tilde{V}_k = 0$ on $\{W_k = 0\}$. Then

$$
\begin{aligned}
(2.3) \qquad V_k &= \left[\frac{1 + \tilde{V}_k}{NX_k}\right]\mathbb{1}_{\{W_k > 0\}} \\
&= \left[\left(\frac{1}{k}\right)\left(\frac{k}{NX_k}\right) + \left(\frac{\tilde{V}_k}{NX_k - k}\right)\left(\frac{NX_k - k}{NX_k}\right)\right]\mathbb{1}_{\{W_k > 0\}}.
\end{aligned}
$$

It follows from (2.2) that

(2.4)     the conditional distribution of $\tilde{V}_k$ given $\mathcal{G}_k$ is Binomial$(NX_k - k, W_k)$

because there are $NX_k - k$ individuals, after the first $k$, with types in $\{1, \ldots, k\}$, and conditional on $\mathcal{G}_k$, each has type $k$ with probability $W_k$. Therefore,

$$
(2.5) \qquad E[V_k|\mathcal{G}_k] = \left[\frac{1}{NX_k} + W_k\left(\frac{NX_k - k}{NX_k}\right)\right]\mathbb{1}_{\{W_k > 0\}}
$$

and

$$
E[V_k|\mathcal{F}_k] = E[E[V_k|\mathcal{G}_k]|\mathcal{F}_k] = r/k.
$$

## 3. Approximating the family-size distribution.

In this section we prove Theorem 1.2, which implies that the distribution $Q_{r,N}$ is a good approximation to the family-size distribution in the Yule process with infinitely many types. To prove this result, we need to show that the $X_k$, which are related to the $V_j$ by (1.1), are close to the $Y_k$, which are likewise related to the $W_j$ by (1.2). We begin by showing that $E[X_k]$ and $E[Y_k]$ are the same.

LEMMA 3.1. *We have* $E[X_k] = E[Y_k] = \prod_{j=k+1}^{N}(1 - \frac{r}{j})$ *for* $1 \le k \le N$.

PROOF. We prove the formula for $E[Y_k]$ by backward induction on $k$. Clearly $E[Y_N] = 1$. Suppose the formula holds for some $k \ge 2$. Then

$$
\begin{aligned}
E[Y_{k-1}] &= E[(1 - W_k)Y_k] = E[(1 - W_k)]E[Y_k] \\
&= \left(1 - \frac{r}{k}\right)\prod_{j=k+1}^{N}\left(1 - \frac{r}{j}\right) = \prod_{j=k}^{N}\left(1 - \frac{r}{j}\right).
\end{aligned}
$$

To get the same formula for $E[X_k]$, first note that $E[X_{k,j}] = 1$ for $1 \le j \le k$. If $n \ge k$, then conditional on $X_{k,n}$, the probability that the $(n+1)$st individual has a type in $\{1, \ldots, k\}$ is $(1-r)X_{k,n}$. Therefore,

$$
E[X_{k,n+1}] = \frac{nE[X_{k,n}] + (1-r)E[X_{k,n}]}{n+1} = \left(1 - \frac{r}{n+1}\right)E[X_{k,n}],
$$

so the formula for $E[X_k]$ follows by induction on $n$.   $\square$

LEMMA 3.2. *We have $(\frac{k}{N})^r e^{-r^2/k} \leq E[X_k] \leq (\frac{k}{N})^r e^{r/k}$ for $1 \leq k \leq N$.*

PROOF. By Lemma 3.1, we have $\log E[X_k] = \sum_{j=k+1}^{N} \log(1 - r/j)$. Note that if $0 \leq x < 1$, then $\log(1-x) = -\sum_{k=1}^{\infty}(x^k/k)$. Summing, we see that if $0 \leq x \leq 1/2$, $-(x + x^2) \leq \log(1-x) \leq -x$. Therefore,

$$
\begin{aligned}
(3.1) \qquad \log E[X_k] &\leq -\sum_{j=k+1}^{N} \frac{r}{j} = \frac{r}{k} - \sum_{j=k}^{N} \frac{r}{j} \\
&\leq \frac{r}{k} - \int_k^N \frac{r}{x}\,dx = \frac{r}{k} + \log\left(\frac{k}{N}\right)^r,
\end{aligned}
$$

$$
\begin{aligned}
(3.2) \qquad \log E[X_k] &\geq -\sum_{j=k+1}^{N} \left(\frac{r}{j} + \frac{r^2}{j^2}\right) \\
&\geq -\int_k^N \frac{r}{x}\,dx - \int_k^N \frac{r^2}{x^2}\,dx \geq \log\left(\frac{k}{N}\right)^r - \frac{r^2}{k}.
\end{aligned}
$$

The result now follows by exponentiating both sides in (3.1) and (3.2). $\square$

LEMMA 3.3. *We have $E[X_k^2(W_k - V_k)^2] \leq r(\frac{1}{N^2} + \frac{2}{N^{1+r}k^{1-r}})$ for $2 \leq k \leq N$.*

PROOF. By (2.3), we have

$$
\begin{aligned}
(3.3) \qquad V_k - W_k &= \left[\left(\frac{1}{k} - W_k\right)\left(\frac{k}{NX_k}\right) \right. \\
&\qquad \left. + \left(\frac{\tilde{V}_k}{NX_k - k} - W_k\right)\left(\frac{NX_k - k}{NX_k}\right)\right]\mathbb{1}_{\{W_k > 0\}}.
\end{aligned}
$$

When we take the conditional expectation given $\mathcal{G}_k$ of the square of the right-hand side of (3.3), the cross-term vanishes because (2.4) implies

$$
E\left[\frac{\tilde{V}_k}{NX_k - k} - W_k \Big| \mathcal{G}_k\right] = 0.
$$

Since $X_k$ and $W_k$ are $\mathcal{G}_k$-measurable, using (2.4) again gives

$$
\begin{aligned}
E[(V_k &- W_k)^2 | \mathcal{G}_k] \\
&= E\left[\left(\frac{1}{k} - W_k\right)^2\left(\frac{k}{NX_k}\right)^2\mathbb{1}_{\{W_k > 0\}} \right. \\
&\qquad \left. + \left(\frac{\tilde{V}_k}{NX_k - k} - W_k\right)^2\left(\frac{NX_k - k}{NX_k}\right)^2\mathbb{1}_{\{W_k > 0\}} \Big| \mathcal{G}_k\right] \\
&= \left(\frac{1}{k} - W_k\right)^2\left(\frac{k}{NX_k}\right)^2\mathbb{1}_{\{W_k > 0\}} + \frac{W_k(1 - W_k)}{NX_k - k}\left(\frac{NX_k - k}{NX_k}\right)^2\mathbb{1}_{\{W_k > 0\}}.
\end{aligned}
$$

Since $W_k$ is independent of $\mathcal{F}_k$ and the conditional distribution of $W_k$ given $W_k > 0$ is Beta$(1, k - 1)$,

$$E[(V_k - W_k)^2 | \mathcal{F}_k] = \frac{k-1}{k^2(k+1)}\left(\frac{k}{NX_k}\right)^2 r + \frac{NX_k - k}{(NX_k)^2}\left(\frac{r}{k} - \frac{2r}{k(k+1)}\right)$$

$$\leq \frac{r}{N^2 X_k^2} + \frac{r}{NkX_k}.$$

Thus, using Lemma 3.2, we get

$$E[X_k^2(W_k - V_k)^2] = E[X_k^2 E[(W_k - V_k)^2 | \mathcal{F}_k]] \leq E\left[\frac{r}{N^2} + \frac{rX_k}{Nk}\right]$$

$$\leq \frac{r}{N^2} + \frac{r}{N^{1+r}k^{1-r}}e^{r/k} \leq r\left(\frac{1}{N^2} + \frac{2}{N^{1+r}k^{1-r}}\right),$$

since for $k \geq 2$, we have $e^{r/k} \leq e^{1/2} \leq 2$.  $\square$

LEMMA 3.4.  *For every real number $a$, we have $E[X_k(X_k - Y_k)(W_k - V_k)(W_k - a)] = 0$.*

PROOF.  Using the fact that $E[W_k - V_k | \mathcal{F}_k] = 0$ and that $W_k$ is $\mathcal{G}_k$-measurable, we have

$$E[(W_k - V_k)(W_k - a)|\mathcal{F}_k] = E[W_k(W_k - V_k)|\mathcal{F}_k]$$

$$= E[E[W_k(W_k - V_k)|\mathcal{G}_k]|\mathcal{F}_k]$$

$$= E[W_k^2 - W_k E[V_k|\mathcal{G}_k]|\mathcal{F}_k].$$

Using (2.5) now, the above equals

$$E\left[W_k^2 - W_k\left(\frac{1}{NX_k} + W_k\left(\frac{NX_k - k}{NX_k}\right)\right)\mathbb{1}_{\{W_k > 0\}}\Big|\mathcal{F}_k\right]$$

$$= \frac{1}{NX_k}(kE[W_k^2|\mathcal{F}_k] - E[W_k|\mathcal{F}_k])$$

$$= \frac{r}{NX_k}\left(\frac{2}{k+1} - \frac{1}{k}\right).$$

It follows that

$$E[X_k(X_k - Y_k)(W_k - V_k)(W_k - a)]$$

$$= E[E[X_k(X_k - Y_k)(W_k - V_k)(W_k - a)|\mathcal{F}_k]]$$

$$= E[X_k(X_k - Y_k)E[(W_k - V_k)(W_k - a)|\mathcal{F}_k]]$$

$$= E\left[(X_k - Y_k)\left(\frac{r}{N}\right)\left(\frac{2}{k+1} - \frac{1}{k}\right)\right] = 0,$$

where the last equality follows from Lemma 3.1.  $\square$

LEMMA 3.5.    *We have $E[(X_k - Y_k)^2] \leq 3/N$ for $1 \leq k \leq N$.*

PROOF.    Suppose $2 \leq k \leq N$. We will bound $E[(X_{k-1} - Y_{k-1})^2]$ in terms of $E[(X_k - Y_k)^2]$. First, note that it follows from (1.1) and (1.2) that

$$
\begin{aligned}
(3.4) \qquad X_{k-1} - Y_{k-1} &= (1 - V_k)X_k - (1 - W_k)Y_k \\
&= X_k(W_k - V_k) + (X_k - Y_k)(1 - W_k).
\end{aligned}
$$

Thus,

$$
\begin{aligned}
(3.5) \qquad E[(X_{k-1} - Y_{k-1})^2] &= E[X_k^2(W_k - V_k)^2] + E[(X_k - Y_k)^2(1 - W_k)^2] \\
&\quad + 2E[X_k(X_k - Y_k)(W_k - V_k)(1 - W_k)].
\end{aligned}
$$

By Lemma 3.4 with $a = 1$, the third term on the right-hand side of (3.5) vanishes. Using Lemma 3.3 and the fact that $E[(X_k - Y_k)^2(1 - W_k)^2] \leq E[(X_k - Y_k)^2]$, we get

$$
E[(X_{k-1} - Y_{k-1})^2] \leq E[(X_k - Y_k)^2] + r\left(\frac{1}{N^2} + \frac{2}{N^{1+r}k^{1-r}}\right).
$$

Since $X_N = Y_N = 1$, it follows that for $1 \leq k \leq N$, we have

$$
\begin{aligned}
(3.6) \qquad E[(X_k - Y_k)^2] &\leq \sum_{j=2}^{N} r\left(\frac{1}{N^2} + \frac{2}{N^{1+r}j^{1-r}}\right) \leq \frac{r}{N} + \frac{2r}{N^{1+r}} \sum_{j=2}^{N} \frac{1}{j^{1-r}} \\
&\leq \frac{1}{N} + \frac{2r}{N^{1+r}} \int_1^N \frac{1}{x^{1-r}}\,dx \leq \frac{1}{N} + \frac{2r}{N^{1+r}}\left(\frac{N^r}{r}\right) = \frac{3}{N},
\end{aligned}
$$

which completes the proof.    $\square$

PROOF OF THEOREM 1.2.    Let $M = \max_{1 \leq k \leq N}|X_k - Y_k|$. Fix $x > 0$. Let $T = \max\{k : |X_k - Y_k| \geq x\}$ if $M \geq x$, and let $T = 0$ otherwise. For $2 \leq k \leq N$, define

$$
(3.7) \qquad \rho_k = X_k(W_k - V_k) - (X_k - Y_k)(W_k - r/k)
$$

so that by (3.4), $X_{k-1} - Y_{k-1} = \rho_k + (X_k - Y_k)(1 - r/k)$. Let

$$
H_k = \begin{cases} X_k - Y_k, & \text{for } k \geq T, \\ X_T - Y_T + \displaystyle\sum_{j=k+1}^{T} \rho_j, & \text{for } k < T. \end{cases}
$$

This definition is chosen so that

$$
(3.8) \qquad H_{k-1} - H_k = \rho_k - (r/k)H_k \mathbb{1}_{\{k > T\}}.
$$

Our first step is to show

$$
(3.9) \qquad P(M \geq x) \leq x^{-2}E[H_1^2].
$$

To establish (3.9), we mimic the proof of Kolmogorov's maximal inequality in [13]. Let $A_k = \{T = k\}$, so the event that $M \geq x$ is the event $\bigcup_{k=1}^{N} A_k$. Then

$$E[H_1^2] \geq \sum_{k=1}^{N} E[H_1^2 \mathbb{1}_{A_k}]$$

$$= \sum_{k=1}^{N} E[(H_k^2 + 2H_k(H_1 - H_k) + (H_1 - H_k)^2) \mathbb{1}_{A_k}]$$

$$\geq \sum_{k=1}^{N} E[H_k^2 \mathbb{1}_{A_k}] + 2 \sum_{k=1}^{N} E[H_k(H_1 - H_k) \mathbb{1}_{A_k}].$$

If $j \leq k$, then

$$\begin{aligned}
(3.10) \quad E[\rho_j | \mathcal{F}_k] &= E[E[\rho_j | \mathcal{F}_j] | \mathcal{F}_k] \\
&= E[X_j E[W_j - V_j | \mathcal{F}_j] - (X_j - Y_j) E[W_j - r/j | \mathcal{F}_j] | \mathcal{F}_k] = 0.
\end{aligned}$$

Therefore,

$$\begin{aligned}
\sum_{k=1}^{N} E[H_k(H_1 - H_k) \mathbb{1}_{A_k}] &= \sum_{k=1}^{N} E[E[H_k(\rho_2 + \cdots + \rho_k) \mathbb{1}_{A_k} | \mathcal{F}_k]] \\
&= \sum_{k=1}^{N} E[H_k \mathbb{1}_{A_k} E[\rho_2 + \cdots + \rho_k | \mathcal{F}_k]] = 0.
\end{aligned}$$

It follows that

$$E[H_1^2] \geq \sum_{k=1}^{N} E[H_k^2 \mathbb{1}_{A_k}] \geq \sum_{k=1}^{N} x^2 P(A_k) = x^2 P(M \geq x),$$

which implies (3.9).

We now obtain a bound on $E[H_1^2]$. Using (3.8) and the fact that the random variable $H_k$ and the event $\{k > T\}$ are $\mathcal{F}_k$-measurable, we have

$$\begin{aligned}
E[H_{k-1}^2 | \mathcal{F}_k] &= E[(\rho_k + H_k(1 - (r/k) \mathbb{1}_{\{k>T\}}))^2 | \mathcal{F}_k] \\
&= E[\rho_k^2 | \mathcal{F}_k] + 2H_k(1 - (r/k) \mathbb{1}_{\{k>T\}}) E[\rho_k | \mathcal{F}_k] \\
&\quad + H_k^2 (1 - (r/k) \mathbb{1}_{\{k>T\}})^2.
\end{aligned}$$

Since $E[\rho_k | \mathcal{F}_k] = 0$ by (3.10), it follows that $E[H_{k-1}^2 | \mathcal{F}_k] \leq E[\rho_k^2 | \mathcal{F}_k] + H_k^2$, and thus $E[H_{k-1}^2] \leq E[\rho_k^2] + E[H_k^2]$. Since $H_N = 0$, we can combine this result with (3.9) to get

$$P(M \geq x) \leq x^{-2} E[H_1^2] \leq x^{-2} \sum_{k=2}^{N} E[\rho_k^2].$$

To bound $E[\rho_k^2]$ we recall the definition in (3.7) and use Lemma 3.5 and the fact that $W_k$ is independent of $X_k$ and $Y_k$ to get

$$E[(X_k - Y_k)^2(W_k - r/k)^2] \leq E[(X_k - Y_k)^2]E\left[W_k^2 - \frac{2r}{k}W_k + \frac{r^2}{k^2}\right]$$

$$\leq \frac{3}{N}\left(\frac{2r}{k(k+1)} - \frac{2r^2}{k^2} + \frac{r^2}{k^2}\right) \leq \frac{6r}{Nk(k+1)}.$$

Combining this result with Lemmas 3.3 and 3.4 with $a = r/k$, we get

$$(3.11) \qquad E[\rho_k^2] \leq r\left(\frac{1}{N^2} + \frac{2}{N^{1+r}k^{1-r}} + \frac{6}{Nk(k+1)}\right).$$

The telescoping sum $\sum_{k=2}^{N} 6r/[Nk(k+1)] \leq 3r/N$, so it follows from (3.6) and (3.11) that

$$P(M \geq x) \leq x^{-2}\sum_{k=2}^{N} E[\rho_k^2] \leq x^{-2}\left(\frac{3}{N} + \frac{3r}{N}\right) \leq \frac{6}{Nx^2}.$$

Thus,

$$E[M] = \int_0^\infty P(M \geq x)\,dx \leq \frac{2}{\sqrt{N}} + \int_{2/\sqrt{N}}^\infty \frac{6}{Nx^2}\,dx = \frac{2}{\sqrt{N}} + \frac{3}{\sqrt{N}} = \frac{5}{\sqrt{N}},$$

which proves the theorem. □

**4. The power law.** In this section we prove Theorem 1.3, which gives the power law for the family-size distribution. Our first lemma gives a bound on the moments of the binomial distribution. Throughout this section, we allow the value of the constant $C$ to change from line to line.

LEMMA 4.1. *Fix $m \geq 1$. There exists a constant $C$ such that for all $n$ and $p$ such that $np \geq 1$, if $X$ has a Binomial$(n, p)$ distribution, then*

$$E\left[\left|\frac{X}{n} - p\right|^m\right] \leq C\left(\frac{p}{n}\right)^{m/2}.$$

PROOF. For now, we assume that $p \leq 1/2$. The proof is based on two bounds for binomial tail probabilities. If $z > 0$, then

$$(4.1) \qquad P\left(\frac{X}{n} - p \leq -z\right) \leq e^{-nz^2/2p},$$

and if $0 < z < 1 - p$, then

$$(4.2) \qquad P\left(\frac{X}{n} - p \geq z\right) \leq e^{-nz^2/2(p+z)}.$$

Equation (4.1) follows from (3.52) on page 121 of [22]. To prove (4.2), we use the fact that if $p < a < 1$, then $P(X/n \geq a) \leq e^{-nH(a)}$, where

$$H(a) = a \log(a/p) + (1-a) \log((1-a)/(1-p)).$$

This is proved, for example, in [5]. We have $H'(a) = \log(a/p) - \log((1-a)/(1-p))$ and $H''(a) = 1/[a(1-a)]$. Since $H(p) = H'(p) = 0$, by Taylor's theorem there exists $z \in [p, a]$ such that $H(a) = \frac{1}{2}H''(z)(a-p)^2$. Note that the function $a \mapsto H''(a)$ is decreasing on $(0, 1/2)$ and increasing on $(1/2, 1)$. Therefore, if $a \leq 1/2$, then $H(a) \geq \frac{1}{2}H''(a)(a-p)^2 \geq \frac{1}{2a}(1-p)^2$ and if $a \geq 1/2$, then $H(a) \geq \frac{1}{2}H''(1/2)(a-p)^2 \geq 2(a-p)^2 \geq \frac{1}{2a}(a-p)^2$. Equation (4.2) follows by substituting $z = a - p$.

Now, using Lemma 5.7 in Chapter 1 of [13], we get

$$
\begin{aligned}
E\left[\left|\frac{X}{n} - p\right|^m\right] = {} & \int_0^p m z^{m-1} P\left(\left|\frac{X}{n} - p\right| > z\right) dz \\
& + \int_p^{1-p} m z^{m-1} P\left(\left|\frac{X}{n} - p\right| > z\right) dz.
\end{aligned}
\tag{4.3}
$$

Using (4.1) and (4.2), then $z \leq p$, and making the substitution $z = y\sqrt{4p/n}$, the first term on the right-hand side is less than or equal to

$$
\begin{aligned}
& \int_0^p m z^{m-1} \left(e^{-nz^2/2p} + e^{-nz^2/2(p+z)}\right) dz \\
& \leq 2m \int_0^p z^{m-1} e^{-nz^2/4p} \, dz \\
& \leq 2m \int_0^\infty \left(\frac{4p}{n}\right)^{m/2} y^{m-1} e^{-y^2} \, dy \\
& \leq C \left(\frac{p}{n}\right)^{m/2}.
\end{aligned}
\tag{4.4}
$$

Likewise, using $z/(p+z) \geq 1/2$ for $z \geq p$ and substituting $z = 4y/n$, the second term on the right-hand side in (4.3) is less than or equal to

$$
\begin{aligned}
& \int_p^{1-p} m z^{m-1} e^{-nz^2/2(p+z)} \, dz \\
& \leq m \int_p^{1-p} z^{m-1} e^{-nz/4} \, dz \\
& \leq m \int_0^\infty \left(\frac{4}{n}\right)^m y^{m-1} e^{-y} \, dy \\
& \leq \frac{C}{n^m}.
\end{aligned}
\tag{4.5}
$$

It follows from (4.3), (4.4) and (4.5) that if $p \leq 1/2$ and $np \geq 1$, then

$$(4.6) \qquad E\left[\left|\frac{X}{n} - p\right|^m\right] \leq C\left(\frac{p}{n}\right)^{m/2} + C\left(\frac{1}{n}\right)^m \leq C\left(\frac{p}{n}\right)^{m/2}.$$

The fact that $np \geq 1$ was used only for the second inequality in (4.6). Therefore, if $p \geq 1/2$ and $np \geq 1$, we can use the first inequality in (4.6) to get

$$E\left[\left|\frac{X}{n} - p\right|^m\right] = E\left[\left|\frac{n-X}{n} - (1-p)\right|^m\right]$$
$$\leq C\left(\frac{1-p}{n}\right)^{m/2} + C\left(\frac{1}{n}\right)^m$$
$$\leq C\left(\frac{p}{n}\right)^{m/2},$$

which completes the proof of the lemma. $\quad\square$

The next lemma bounds the moments of $X_k$. Recall that $X_k = X_{k,N}$ is the fraction of the first $N$ individuals with one of the first $k$ types.

LEMMA 4.2. *Fix a real number $m \geq 1$. Then there is a positive constant $C$ such that for all $k \geq 1$,*

$$E[X_k^m] \leq \left(\frac{k}{N}\right)^{mr}\left(1 + \frac{C}{k}\right).$$

PROOF. Let $M_{k,l} = \sum_{j=1}^k R_{j,l}$ be the number of individuals at time $T_l$ with types in $\{1,\ldots,k\}$. Note that $M_{k,N} = NX_k$. Conditional on $M_{k,l}$, the probability that the $(l+1)$st individual born has a type in $\{1,\ldots,k\}$ is $(1-r)M_{k,l}/l$. Therefore,

$$E[M_{k,l+1}^m | M_{k,l}] = M_{k,l}^m + (1-r)\left(\frac{M_{k,l}}{l}\right)[(M_{k,l}+1)^m - M_{k,l}^m].$$

Since $b^m - a^m = \int_a^b mx^{m-1}\,dx \leq mb^{m-1}(b-a)$ for $0 \leq a \leq b$, the above is less than or equal to

$$M_{k,l}^m + (1-r)\left(\frac{M_{k,l}}{l}\right)m(M_{k,l}+1)^{m-1}$$
$$= M_{k,l}^m + \frac{(1-r)m}{l}M_{k,l}^m + \frac{(1-r)m}{l}[(M_{k,l}+1)^{m-1} - M_{k,l}^{m-1}]M_{k,l}.$$

Using the integration inequality again this is less than or equal to

$$M_{k,l}^m\left(1 + \frac{(1-r)m}{l}\right) + \frac{(1-r)m}{l}[(m-1)(M_{k,l}+1)^{m-2}]M_{k,l}$$
$$\leq M_{k,l}^m\left(1 + \frac{(1-r)m}{l}\right) + \frac{m(m-1)}{l}(M_{k,l}+1)^{m-1}.$$

Since $M_{k,l} \geq 1$, we have

$$(4.7) \qquad E[M_{k,l+1}^m | M_{k,l}] \leq M_{k,l}^m \left(1 + \frac{(1-r)m}{l}\right) + \frac{C}{l} M_{k,l}^{m-1}.$$

We now establish the lemma for integer values of $m$ by induction. When $m = 1$, the result is an immediate consequence of Lemma 3.2 and the inequality $e^{r/k} \leq 1 + C/k$. Suppose the result holds for $m - 1$. Then, since $M_{k,l} = lX_{k,l}$, we have

$$E[M_{k,l}^{m-1}] = l^{m-1} E\left[\left(\frac{M_{k,l}}{l}\right)^{m-1}\right]$$

$$\leq l^{m-1} \left(\frac{k}{l}\right)^{(m-1)r} \left(1 + \frac{C}{k}\right)$$

$$\leq C k^{(m-1)r} l^{(m-1)(1-r)}.$$

Therefore, taking expectations of both sides in (4.7), we get

$$E[M_{k,l+1}^m] \leq \left(1 + \frac{(1-r)m}{l}\right) E[M_{k,l}^m] + C k^{(m-1)r} l^{(m-1)(1-r)-1}.$$

Since $M_{k,k} = k$, iterating the last result shows that $E[X_k^m] = E[M_{k,N}]/N^m$ is at most

$$\frac{1}{N^m}\left[k^m \prod_{j=k}^{N-1} \left(1 + \frac{(1-r)m}{j}\right)\right.$$

$$\left. + \sum_{l=k}^{N-1} C k^{(m-1)r} l^{(m-1)(1-r)-1} \left(\prod_{j=l+1}^{N-1} \left(1 + \frac{(1-r)m}{j}\right)\right)\right].$$

Since $1 + x \leq e^x$ for $x > 0$, we have

$$\prod_{j=k}^{N-1} \left(1 + \frac{(1-r)m}{j}\right) \leq \exp\left(\sum_{j=k}^{N-1} \frac{(1-r)m}{j}\right)$$

$$\leq \exp\left((1-r)m\left(\frac{1}{k} + \int_k^N x^{-1}\,dx\right)\right)$$

$$= \exp\left(\frac{(1-r)m}{k} + (1-r)m \log\left(\frac{N}{k}\right)\right)$$

$$\leq \left(\frac{N}{k}\right)^{(1-r)m} \left(1 + \frac{C}{k}\right).$$

Thus,

$$E[X_k^m] \leq \frac{1}{N^m}\left[k^m \left(\frac{N}{k}\right)^{(1-r)m} \left(1 + \frac{C}{k}\right)\right.$$

$$+ C \sum_{l=k}^{N-1} k^{(m-1)r} l^{(m-1)(1-r)-1} \left( \frac{N}{l} \right)^{(1-r)m} \Bigg]$$

$$= \left( \frac{k}{N} \right)^{mr} \left[ 1 + \frac{C}{k} + Ck^{-r} \sum_{l=k}^{N-1} l^{-2+r} \right] \le \left( \frac{k}{N} \right)^{mr} \left[ 1 + \frac{C}{k} \right].$$

The result for integer values of $m$ follows by induction.

Now suppose $n < m < n + 1$, where $n$ is a positive integer. Let

$$p = (n - m + 1)^{-1}$$

and let

$$q = (m - n)^{-1}.$$

Note that $p^{-1} + q^{-1} = 1$ and $n/p + (n+1)/q = m$. By Hölder's inequality,

$$\begin{aligned}
E[X_k^m] &= E[X_k^{n/p} X_k^{(n+1)/q}] \\
&\le E[X_k^n]^{m-n+1} E[X_k^{n+1}]^{m-n} \\
&\le \left( \frac{k}{N} \right)^{mr} \left( 1 + \frac{C}{k} \right),
\end{aligned}$$

so the lemma is true for all real numbers $m \ge 1$. $\quad\square$

To prove Theorem 1.3, we will approximate the family sizes $NV_k X_k$ by $NW_k(k/N)^r$. To use this approximation, we will need a bound on the probability that the difference between these two quantities is large. Note that

$$(4.8) \qquad V_k X_k - W_k \left( \frac{k}{N} \right)^r = X_k(V_k - W_k) + W_k \left( X_k - \left( \frac{k}{N} \right)^r \right).$$

The next two lemmas deal separately with the two terms on the right-hand side of (4.8).

LEMMA 4.3. *There is a positive constant $C$ so that for all $\delta > 0$*

$$\sum_{k=1}^N P \left( \left| W_k \left( X_k - \left( \frac{k}{N} \right)^r \right) \right| > \frac{\delta S}{2N} \right) \le C \left( \frac{N^{1-r}}{\delta S} \right)^{2/(3-2r)}.$$

PROOF. Conditioning on $\mathcal{F}_k$ and noting that $W_k$ is independent of $\mathcal{F}_k$ gives

$$E \left[ W_k^2 \left( X_k - \left( \frac{k}{N} \right)^r \right)^2 \right] = E[W_k^2] E \left[ \left( X_k - \left( \frac{k}{N} \right)^r \right)^2 \right].$$

If $k \geq 2$, Lemmas 3.2 and 4.2 give that the above is equal to

$$\frac{2r}{k(k+1)}\left(E[X_k^2] - 2E[X_k]\left(\frac{k}{N}\right)^r + \left(\frac{k}{N}\right)^{2r}\right)$$

$$\leq \frac{2r}{k^2}\left[\left(\frac{k}{N}\right)^{2r}\left(1 + \frac{C}{k}\right) - 2\left(\frac{k}{N}\right)^{2r}e^{-r^2/k} + \left(\frac{k}{N}\right)^{2r}\right]$$

$$\leq \frac{2r}{k^2}\left(\frac{k}{N}\right)^{2r}\left[1 + \frac{C}{k} - 2\left(1 - \frac{r^2}{k}\right) + 1\right]$$

$$\leq \frac{C}{N^{2r}k^{3-2r}}.$$

Fix a positive integer $L$. Using a trivial inequality for $k \leq L$ and Chebyshev's inequality,

$$\sum_{k=1}^{N} P\left(\left|W_k\left(X_k - \left(\frac{k}{N}\right)^r\right)\right| > \frac{\delta S}{2N}\right) \leq L + \sum_{k=L+1}^{N} \frac{C}{N^{2r}k^{3-2r}}\left(\frac{\delta S}{2N}\right)^{-2}$$

$$(4.9) \qquad\qquad\qquad \leq L + \frac{CN^{2-2r}}{(\delta S)^2}\int_L^\infty \frac{1}{x^{3-2r}}\,dx$$

$$= L + \frac{CN^{2-2r}L^{-(2-2r)}}{(2-2r)(\delta S)^2}.$$

If $L = (N^{1-r}/(\delta S))^{2/(3-2r)}$, then the right-hand side of (4.9) is bounded by

$$\left(\frac{N^{1-r}}{\delta S}\right)^{2/(3-2r)} + C\left(\frac{N^{1-r}}{\delta S}\right)^{2-(2-2r)2/(3-2r)} \leq C\left(\frac{N^{1-r}}{\delta S}\right)^{2/(3-2r)},$$

as claimed. □

Lemma 4.4. *There is a constant $C$ so that for all $\delta > 0$, we have*

$$(4.10) \qquad \sum_{k=1}^{N} P\left(|X_k(V_k - W_k)| > \frac{\delta S}{2N}\right) \leq 1 + \frac{CN}{(\delta S)^{3/2(1-r)}}.$$

Proof. Recall from Section 2 that $W_k = \xi_k \tilde{W}_k$, where $\xi_k = \mathbb{1}_{\{W_k > 0\}}$ has a Bernoulli($r$) distribution and $\tilde{W}_k$ has a Beta($1, k-1$) distribution and is independent of $\xi_k$. Also recall that $\tilde{V}_k = (NV_kX_k - 1)\mathbb{1}_{\{W_k > 0\}}$ is a random variable such that the conditional distribution of $\tilde{V}_k$ given $\mathcal{G}_k$ is Binomial($NX_k - k, \tilde{W}_k$). Using (2.3), we see that for all $k \geq 2$ we have

$$P\left(|X_k(V_k - W_k)| > \frac{\delta S}{2N}\right)$$

$$= P\left(|NX_k(V_k - W_k)|\mathbb{1}_{\{W_k > 0\}} > \frac{\delta S}{2}\right)$$

(4.11)
$$= P\left(|(1 - k\tilde{W}_k) + (\tilde{V}_k - \tilde{W}_k(NX_k - k))|\mathbb{1}_{\{W_k > 0\}} > \frac{\delta S}{2}\right)$$

$$\leq P\left(|1 - k\tilde{W}_k| > \frac{\delta S}{4}\right) + P\left(|\tilde{V}_k - \tilde{W}_k(NX_k - k)| > \frac{\delta S}{4}\right).$$

Let $m = 3/2(1-r)$. The reason for this choice will become clear in (4.15). Until then the reader should keep in mind that $m$ is a fixed real number. Since $\Gamma(x + 1) = x\Gamma(x)$ for all real $x$, we have $\Gamma(k)/\Gamma(m + k) \leq Ck^{-m}$ for some constant $C$. Therefore,

(4.12)
$$E[\tilde{W}_k^m] = \frac{\Gamma(k)\Gamma(m+1)}{\Gamma(m+k)} \leq Ck^{-m},$$

so using $(a + b)^m \leq 2^m(a^m + b^m)$ for $a, b \geq 0$, we have

$$E[(1 + k\tilde{W}_k)^m] \leq 2^m(1 + E[(k\tilde{W}_k)^m]) \leq C.$$

Therefore, by Markov's inequality, if $k \geq 2$, then

$$P\left(|1 - k\tilde{W}_k| > \frac{\delta S}{4}\right) \leq P\left(|1 + k\tilde{W}_k| > \frac{\delta S}{4}\right)$$

(4.13)
$$\leq \left(\frac{\delta S}{4}\right)^{-m} E[(1 + k\tilde{W}_k)^m]$$

$$\leq \frac{C}{(\delta S)^m},$$

which bounds the first term on the right-hand side of (4.11).

Because of the restriction $np \geq 1$ in Lemma 4.1, we must split the second term in (4.11) into two pieces, depending on the value of $\tilde{W}_k(NX_k - k)$. Let $V'_k$ be a random variable such that, conditional on $\mathcal{G}_k$, the distribution of $V'_k$ is Binomial$(NX_k - k, 1/(NX_k - k))$. We set $V'_k = 0$ if $NX_k - k = 0$. Note that when $\tilde{W}_k(NX_k - k) < 1$, the conditional distribution of $V'_k$ given $\mathcal{G}_k$ stochastically dominates the conditional distribution of $\tilde{V}_k$ given $\mathcal{G}_k$. By Lemma 4.1, $E[|V'_k - 1|^m|\mathcal{G}_k] \leq C$. Note also that $|\tilde{V}_k - \tilde{W}_k(NX_k - k)| = 0$ on the event $\{NX_k - k = 0\}$. Therefore, if $k \geq 2$, then

$$P\left(|\tilde{V}_k - \tilde{W}_k(NX_k - k)|\mathbb{1}_{\{\tilde{W}_k(NX_k - k) < 1\}} > \frac{\delta S}{4}\right)$$

$$= E\left[E\left[P\left(|\tilde{V}_k - \tilde{W}_k(NX_k - k)|\mathbb{1}_{\{\tilde{W}_k(NX_k - k) < 1\}} > \frac{\delta S}{4}\right)\Big|\mathcal{G}_k\right]\right]$$

(4.14)
$$\leq E\left[E\left[P\left(|V'_k - 1| + 1 > \frac{\delta S}{4}\right)\Big|\mathcal{G}_k\right]\right]$$

$$\leq \left(\frac{\delta S}{4}\right)^{-m} E[E[(|V'_k - 1| + 1)^m|\mathcal{G}_k]]$$

$$\leq \left(\frac{\delta S}{4}\right)^{-m} 2^m (C+1) \leq \frac{C}{(\delta S)^m}.$$

By Lemma 4.1, we get, for $k \geq 2$,

$$P\left(|\tilde{V}_k - \tilde{W}_k(NX_k - k)|\mathbb{1}_{\{\tilde{W}_k(NX_k-k)\geq 1\}} > \frac{\delta S}{4}\right)$$

$$= E\left[E\left[P\left(\left|\frac{\tilde{V}_k}{NX_k - k} - \tilde{W}_k\right|\mathbb{1}_{\{\tilde{W}_k(NX_k-k)\geq 1\}} > \frac{\delta S}{4(NX_k - k)}\right)\Big|\mathcal{G}_k\right]\right]$$

$$\leq E\left[\left(\frac{\delta S}{4(NX_k - k)}\right)^{-m} E\left[\left|\frac{\tilde{V}_k}{NX_k - k} - \tilde{W}_k\right|^m \mathbb{1}_{\{\tilde{W}_k(NX_k-k)\geq 1\}}\Big|\mathcal{G}_k\right]\right]$$

$$\leq E\left[\left(\frac{4(NX_k - k)}{\delta S}\right)^m C\left(\frac{\tilde{W}_k}{NX_k - k}\right)^{m/2}\right]$$

$$\leq \frac{C}{(\delta S)^m} E[\tilde{W}_k^{m/2}(NX_k - k)^{m/2}].$$

By conditioning on $\mathcal{F}_k$ and noting that $W_k$ is independent of this $\sigma$-field, we see that this is at most

$$\frac{C}{(\delta S)^m} E[\tilde{W}_k^{m/2}] E[(NX_k)^{m/2}] \leq \frac{C}{(\delta S)^m}\frac{C}{k^{m/2}} N^{m/2}\left(\frac{k}{N}\right)^{mr/2}$$

by (4.12) and Lemma 4.2. Recalling that $m = 3/2(1 - r)$, the above is at most

$$(4.15) \qquad \frac{C}{(\delta S)^m}\left(\frac{N}{k}\right)^{m(1-r)/2} = \frac{C}{(\delta S)^m}\left(\frac{N}{k}\right)^{3/4}.$$

Note that

$$\sum_{k=2}^{N}\left(\frac{N}{k}\right)^{3/4} = N^{3/4}\sum_{k=2}^{N} k^{-3/4} \leq CN^{3/4}N^{1/4} = CN.$$

Combining this fact with (4.11), (4.13), (4.14) and (4.15), which hold for $k \geq 2$, we get

$$\sum_{k=1}^{N} P\left(|X_k(V_k - W_k)| > \frac{\delta S}{2N}\right) \leq 1 + \frac{CN}{(\delta S)^m},$$

which completes the proof.   $\square$

Now that we have shown that $V_k X_k$ and $W_k(k/N)^r$ are close with high probability, the next step is to calculate the probability that $W_k(k/N)^r$ is large. The next two lemmas provide upper and lower bounds on the probability that $W_k(k/N)^r$ is large. Recall from (1.3) that

$$g(S) = r\Gamma\left(\frac{2-r}{1-r}\right)NS^{-1/(1-r)}.$$

LEMMA 4.5. *There is a constant $C$ so that for $0 < \delta < 1/2$ and $S \leq \frac{1}{2}N^{1-r}$,*

$$\sum_{k=1}^{N} P\left(W_k\left(\frac{k}{N}\right)^r \geq \frac{(1-\delta)S}{N}\right) \leq C + g(S)(1 + C\delta).$$

PROOF. The $C$ on the left takes care of the term $k = 1$. Since the conditional distribution of $W_k$ given $W_k > 0$ is Beta$(1, k-1)$, we have, for all $k \geq 2$ and $a \in (0, 1)$,

$$(4.16) \qquad P(W_k \geq a) = r \int_a^1 (k-1)(1-x)^{k-2} \, dx = r(1-a)^{k-1}.$$

Using the facts that $(1 - a/x)^x \leq e^{-a}$ if $0 \leq a \leq x$ and $1/(1-x) \leq 1 + 2x$ if $0 \leq x \leq 1/2$, we have, for $k \geq 2$,

$$(4.17) \quad \begin{aligned} P\left(W_k \geq \frac{S(1-\delta)}{k^r N^{1-r}}\right) &= r\left(1 - \frac{S(1-\delta)k^{1-r}}{kN^{1-r}}\right)^k\left(1 - \frac{S(1-\delta)}{k^r N^{1-r}}\right)^{-1} \\ &\leq re^{-S(1-\delta)(k/N)^{1-r}}\left(1 + \frac{2S}{k^r N^{1-r}}\right). \end{aligned}$$

Note that

$$(4.18) \quad \begin{aligned} \sum_{k=2}^{N} &e^{-S(1-\delta)(k/N)^{1-r}}\left(1 + \frac{2S}{k^r N^{1-r}}\right) \\ &\leq \int_0^{\infty} e^{-S(1-\delta)(x/N)^{1-r}}\left(1 + \frac{2S}{x^r N^{1-r}}\right) dx. \end{aligned}$$

Letting $y = S(1-\delta)(x/N)^{1-r}$, which means that $x = y^{1/(1-r)}M$ where $M = N(S(1-\delta))^{-1/(1-r)}$ and $dx = y^{1/(1-r)-1}M/(1-r) \, dy$, we have

$$(4.19) \quad \begin{aligned} \int_0^{\infty} e^{-S(1-\delta)(x/N)^{1-r}} \, dx &= \int_0^{\infty} e^{-y} y^{1/(1-r)-1}\left(\frac{M}{1-r}\right) dy \\ &= \Gamma\left(\frac{1}{1-r}\right)\frac{M}{1-r} \\ &= \Gamma\left(\frac{2-r}{1-r}\right)N(S(1-\delta))^{-1/(1-r)}. \end{aligned}$$

The same change of variables gives

$$\begin{aligned} \int_0^{\infty} &e^{-S(1-\delta)(x/N)^{1-r}}\left(\frac{2S}{x^r N^{1-r}}\right) dx \\ &= \frac{2S}{N^{1-r}} \int_0^{\infty} e^{-y}(y^{1/(1-r)}M)^{-r} y^{1/(1-r)-1}\left(\frac{M}{1-r}\right) dy \end{aligned}$$

(4.20)
$$\begin{aligned}
&= \frac{2SM^{1-r}}{(1-r)N^{1-r}} \int_0^\infty e^{-y}\,dy \\
&= \frac{2}{(1-r)(1-\delta)} \leq C.
\end{aligned}$$

Because $(1-\delta)^{-1/(1-r)} \leq 1 + C\delta$, the claim follows from (4.17)–(4.20).   □

Lemma 4.6.  *There is a constant $C$ so that for $0 < \delta < 1/2$ and $S \leq \frac{1}{3}N^{1-r}$,*

$$\sum_{k=1}^N P\left(W_k\left(\frac{k}{N}\right)^r \geq \frac{(1+\delta)S}{N}\right) \geq -C + g(S)(1 - C\delta - Ce^{-S/2}).$$

Proof.  Recall from the beginning of the proof of Lemma 3.2 that if $0 \leq x \leq 1/2$, then $\log(1-x) \geq -(x+x^2)$. It follows that if $0 \leq a/y \leq 1/2$ and $a \geq 0$, then $y\log(1-a/y) \geq -a - a^2/y$, and so

$$\left(1 - \frac{a}{y}\right)^y \geq e^{-a}e^{-a^2/y} \geq e^{-a}\left(1 - \frac{a^2}{y}\right).$$

Therefore, if $k \geq 2$, then

$$\begin{aligned}
&P\left(W_k\left(\frac{k}{N}\right)^r \geq \frac{(1+\delta)S}{N}\right) \\
&\quad = r\left(1 - \frac{S(1+\delta)}{k^r N^{1-r}}\right)^{k-1} \\
&\quad \geq r\left(1 - \frac{S(1+\delta)k^{1-r}}{kN^{1-r}}\right)^k \\
&\quad \geq re^{-S(1+\delta)(k/N)^{1-r}}\left(1 - \frac{S^2(1+\delta)^2 k^{1-2r}}{N^{2-2r}}\right).
\end{aligned}$$

(4.21)

We have

$$\begin{aligned}
&\sum_{k=2}^N e^{-S(1+\delta)(k/N)^{1-r}} \\
&\quad \geq \left(\int_0^\infty e^{-S(1+\delta)(x/N)^{1-r}}\,dx\right) \\
&\quad \quad - 2 - \int_N^\infty e^{-S(1+\delta)(x/N)^{1-r}}\,dx.
\end{aligned}$$

(4.22)

It follows from (4.19) with $\delta$ replaced by $-\delta$ and $M = N(S(1+\delta))^{-1/(1-r)}$ that

$$(4.23) \qquad \int_0^\infty e^{-S(1+\delta)(x/N)^{1-r}}\,dx = \Gamma\left(\frac{2-r}{1-r}\right)M$$

$$\geq \Gamma\left(\frac{2-r}{1-r}\right)NS^{-1/(1-r)}(1-C\delta).$$

To estimate the second term in (4.21) we note that

$$\sum_{k=2}^N e^{-S(1+\delta)(k/N)^{1-r}}\left(\frac{S^2(1+\delta)^2 k^{1-2r}}{N^{2-2r}}\right)$$

$$\leq \frac{CS^2}{N^{2-2r}}\int_0^\infty e^{-S(1+\delta)(x/N)^{1-r}}x^{1-2r}\,dx.$$

Making the change of variables $x = y^{1/(1-r)}M$ and reasoning as in (4.20), we see that this equals

$$(4.24) \qquad \frac{CS^2}{N^{2-2r}}\int_0^\infty e^{-y}(y^{1/(1-r)}M)^{1-2r}y^{1/(1-r)-1}\left(\frac{M}{1-r}\right)dy$$

$$= \frac{C}{(1+\delta)^2}\int_0^\infty e^{-y}y\,dy \leq C.$$

For all real numbers $b$, there is a constant $C$ such that $x^b e^{-x/2} \leq C$ for all $x > 1$. Using this fact and our favorite change of variables, we get

$$(4.25) \qquad \int_N^\infty e^{-S(1+\delta)(x/N)^{1-r}}\,dx = \int_{S(1+\delta)}^\infty e^{-y}y^{1/(1-r)-1}\left(\frac{M}{1-r}\right)dy$$

$$\leq \frac{M}{1-r}\int_S^\infty e^{-y}y^{1/(1-r)-1}\,dy$$

$$\leq \frac{CM}{1-r}\int_S^\infty e^{-y/2}\,dy$$

$$\leq CNS^{-1/(1-r)}e^{-S/2}.$$

The lemma now follows by combining (4.21)–(4.25). $\square$

Proof of Theorem 1.3. Let $A_k$ be the event that $\{NX_kV_k \geq S\}$, so $F_{S,N} = \sum_{k=1}^N \mathbb{1}_{A_k}$. First, note that if the theorem is true for $S = \frac{1}{3}N^{1-r}$, then we know that $E[F_{S,N}] \leq C$ for all $S > \frac{1}{3}N^{1-r}$, which implies the assertion in the theorem. Therefore, it suffices to prove the result for $S \leq \frac{1}{3}N^{1-r}$, in which case the conclusions of Lemmas 4.5 and 4.6 will hold as long as we choose $\delta < 1/2$. Let $A_k^- = \{NW_k(k/N)^r \geq (1-\delta)S\}$ and let $A_k^+ = \{NW_k(k/N)^r \geq$

$(1 + \delta)S\}$. Let $F_S^- = \sum_{k=1}^N \mathbb{1}_{A_k^-}$ and $F_S^+ = \sum_{k=1}^N \mathbb{1}_{A_k^+}$. Writing $F_S$ for $F_{S,N}$, we have

$$(4.26) \quad |F_S - g(S)| \le |F_S - F_S^-| + |F_S^- - E[F_S^-]| + |E[F_S^-] - g(S)|.$$

To prove the theorem, we will bound the expectations of the three terms on the right-hand side of (4.26).

Note that $A_k^+ \subset A_k^-$ for all $k$ and $A_k \triangle A_k^- \subset (A_k^- \setminus A_k^+) \cup \{|V_k X_k - W_k(k/N)^r| > \delta S/N\}$. By Lemmas 4.5 and 4.6, we have

$$(4.27) \qquad \sum_{k=1}^N P(A_k^- \setminus A_k^+) \le C + Cg(S)\{\delta + e^{-S/2}\}.$$

By (4.8) and Lemmas 4.3 and 4.4, we have

$$\sum_{k=1}^N P\left(\left|V_k X_k - W_k\left(\frac{k}{N}\right)^r\right| > \frac{\delta S}{N}\right)$$

$$(4.28) \qquad \le 1 + C\left(\frac{N^{1-r}}{\delta S}\right)^{2/(3-2r)} + C\left(\frac{N}{(\delta S)^{3/2(1-r)}}\right)$$

$$\le 1 + Cg(S)\left\{\frac{(NS^{-1/(1-r)})^{-1/(3-2r)}}{\delta^{2/(3-2r)}} + \frac{S^{-1/2(1-r)}}{\delta^{3/2(1-r)}}\right\}.$$

Combining the last two results, we get

$$(4.29) \qquad E[|F_S - F_S^-|] \le C + Cg(S)(D_1 + D_2),$$

where $D_1$ and $D_2$ are the terms in braces in (4.27) and ( 4.28). To bound the second term of (4.26), we use Jensen's inequality and the fact that the $A_k^-$ are independent to get

$$E[|F_S^- - E[F_S^-]|] \le E[(F_S^- - E[F_S^-])^2]^{1/2}$$

$$= \mathrm{Var}(F_S^-)^{1/2} = \left[\sum_{k=1}^N \mathrm{Var}(\mathbb{1}_{A_k^-})\right]^{1/2}$$

$$(4.30) \qquad \le \left[\sum_{k=1}^N P(A_k^-)\right]^{1/2}$$

$$\le Cg(S)^{1/2} \le Cg(S)(NS^{-1/(1-r)})^{-1/2}.$$

Furthermore, note that since Lemma 4.5 gives an upper bound for $E[F_S^-]$ that is greater than $g(S)$ and Lemma 4.6 gives a lower bound for $E[F_S^+]$ that is smaller than $g(S)$, the difference $|E[F_S^-] - g(S)|$ is less than or equal to the difference between these two bounds, which itself was bounded in (4.27). Combining this observation with (4.26), (4.29) and (4.30), we see that

$$E[|F_S - g(S)|] \leq C + Cg(S)(D_1 + D_2 + (NS^{-1/(1-r)})^{-1/2}).$$

To prove the theorem, we need to show that each part of $D_1 + D_2$ is bounded by $S^{-a}$ or $(NS^{-1/(1-r)})^{-b}$ for some positive constants $a$ and $b$. Letting $R = NS^{-1/(1-r)}$ to simplify notation, it is enough to bound

$$\delta + \frac{R^{-1/(3-2r)}}{\delta^{2/(3-2r)}} + \frac{S^{-1/2(1-r)}}{\delta^{3/2(1-r)}}.$$

To do this, we let $\delta = A(S^{-c} + R^{-d})$, where $3c < 1$ and $2d < 1$, and choose $A$ to ensure that $\delta < 1/2$. To optimize the bound we set $(1 - 3c)/2(1 - r) = c$ and $(1 - 2d)/(3 - 2r) = d$. Solving gives $c = 1/(5 - 2r)$ and $d = 1/(5 - 2r)$. $\square$

**5. Sizes of the largest families.** In this section we study the largest families, whose sizes are $O(N^{1-r})$, and we prove Propositions 1.4 and 1.5. The key to our arguments is the following well-known result about Yule processes, which is discussed in Chapter III of [7]. Suppose $(X(t), t \geq 0)$ is a Yule process started with one individual at time zero in which each individual splits into two at rate $\lambda$. Then, there exists a random variable $W$ such that

$$\lim_{t \to \infty} e^{-\lambda t} X(t) = W \qquad \text{a.s.}$$

and $W$ has an exponential distribution with mean 1. A consequence of this fact is that if $X_1(t), \ldots, X_k(t)$ are $k$ independent Yule processes, each started with one individual at time zero, then

$$(5.1) \qquad \lim_{t \to \infty} \frac{X_1(t)}{X_1(t) + \cdots + X_k(t)} = B \qquad \text{a.s.},$$

where $B$ has the Beta$(1, k - 1)$ distribution.

PROOF OF PROPOSITION 1.4. The $k = 1$ case was proved by Angerer [2], so we may fix $k \geq 2$. Let $\mathcal{I}_k$ denote the $k$th individual to enter the population. Let $D_{k,N}$ be the number of descendants of $\mathcal{I}_k$ in the population at time $T_N$, when the total population size reaches $N$, and let $G_{k,N}$ be the number of those descendants having the same type as $\mathcal{I}_k$. It follows from (5.1) that

$$\lim_{N \to \infty} \frac{D_{k,N}}{N} = B_k \qquad \text{a.s.},$$

where $B_k$ has the Beta$(1, k - 1)$ distribution. Also, by the same argument as in the $k = 1$ case, we have

$$\lim_{N \to \infty} \frac{G_{k,N}}{(D_{k,N})^{1-r}} = M_k \qquad \text{a.s.},$$

where $M_k$ has the Mittag–Leffler distribution with parameter $1 - r$. Moreover, since the descendants of $\mathcal{I}_k$ form a Yule process and mutations are neutral, $M_k$ and $B_k$ are independent.

Recall that $R_{k,N}$ is the number of type-$k$ individuals in the population at time $T_N$. On the event that the $k$th individual born is a mutant, we have $R_{k,N} = G_{k,N}$. Therefore

$$Z_k = \lim_{N \to \infty} \frac{R_{k,N}}{N^{1-r}} = \lim_{N \to \infty} \frac{G_{k,N}}{(D_{k,N})^{1-r}} \left( \frac{D_{k,N}}{N} \right)^{1-r} = M_k B_k^{1-r}$$

almost surely on the event that the $k$th individual born is a mutant. Proposition 1.4 follows because $M_k$ and $B_k$ are independent of the event that the $k$th individual is a mutant.  □

It remains to prove Proposition 1.5. We will need the following lemma.

LEMMA 5.1.   *Given $\varepsilon > 0$, there exists a positive integer $L$ such that for sufficiently large $N$,*

$$\sum_{k=L}^{N} P(R_{k,N} \geq N^{1-r}) < \varepsilon. \tag{5.2}$$

PROOF.   We have $R_{k,N} = N V_k X_k$, so $P(R_{k,N} \geq N^{1-r}) = P(V_k X_k \geq N^{-r})$. From ( 4.9) with $S = N^{1-r}$ and $\delta = 1/2$, we get

$$\sum_{k=L}^{N} P\left( \left| W_k \left( X_k - \left( \frac{k}{N} \right)^r \right) \right| > \frac{N^{-r}}{4} \right) \leq \frac{2C}{1-r}(L-1)^{-(2-2r)}, \tag{5.3}$$

which is less than $\varepsilon/3$ for sufficiently large $L$. By Lemma 4.4, again with $S = N^{1-r}$ and $\delta = 1/2$,

$$\sum_{k=L}^{N} P\left( |X_k(V_k - W_k)| > \frac{N^{-r}}{4} \right) < \frac{\varepsilon}{3} \tag{5.4}$$

for sufficiently large $N$ as long as $L \geq 2$, because the 1 on the right-hand side of (4.10) comes from the $k = 1$ term. Finally, (4.16) implies that for $L \geq 2$,

$$\sum_{k=L}^{N} P\left( W_k \left( \frac{k}{N} \right)^r > \frac{N^{-r}}{2} \right) = \sum_{k=L}^{N} P\left( W_k > \frac{1}{2k^r} \right) = \sum_{k=L}^{N} r \left( 1 - \frac{1}{2k^r} \right)^{k-1},$$
$$\tag{5.5}$$

which is also at most $\varepsilon/3$ for sufficiently large $L$. The lemma follows from (5.3), (5.4) and (5.5).  □

We now review some facts about the Mittag–Leffler distribution. Let $X$ be a stable random variable satisfying $E[e^{-\lambda X}] = e^{-\lambda^\alpha}$, where $0 < \alpha < 1$.

Then, it is well known (see [30], Section 0.5) that $X$ is nonnegative and has density

$$(5.6) \qquad f_\alpha(x) = \frac{1}{\pi} \sum_{k=0}^{\infty} \frac{(-1)^{k-1}}{k!} \sin(\pi\alpha k) \frac{\Gamma(\alpha k + 1)}{x^{\alpha k+1}}, \qquad x > 0.$$

It follows from [36] that if $A_1 = \alpha^{1/(2(1-\alpha))} (\cos \frac{\pi\alpha}{2})^{-1/(2(1-\alpha))} [2\pi(1-\alpha)]^{-1/2}$ and $A_2 = (1-\alpha)\alpha^{\alpha/(1-\alpha)} (\cos \frac{\pi\alpha}{2})^{-1/(1-\alpha)}$, then

$$f_\alpha(x) \sim A_1 x^{-1-\alpha/(2(1-\alpha))} \exp(-A_2 x^{-\alpha/(1-\alpha)}),$$

where "$\sim$" means that the ratio of the two sides tends to 1 as $x \downarrow 0$. The Mittag–Leffler distribution with parameter $\alpha \in (0,1)$ is the distribution of $Y = X^{-\alpha}$. Therefore, if $g_\alpha$ denotes the density of $Y$, a change of variables gives

$$(5.7) \qquad g_\alpha(x) = \frac{f_\alpha(x^{-1/\alpha})}{\alpha x^{1+1/\alpha}} \sim \frac{A_1}{\alpha} x^{1/(2(1-\alpha))-1} \exp(-A_2 x^{1/(1-\alpha)}),$$

where "$\sim$" means that the ratio of the two sides tends to 1 as $x \to \infty$. In the following proof, $C$ is a positive constant whose value may change from line to line.

PROOF OF PROPOSITION 1.5. Let $g$ be the density of $M$, which has the Mittag–Leffler distribution with parameter $1 - r$. By (5.7), there exists a constant $C$ such that $g(x) \leq C x^{1/2r-1} e^{-A_2 x^{1/r}}$ for all $x \geq 1$. Therefore, if $x \geq 1$, then, making the substitution $y = A_2 z^{1/r}$, we get

$$\begin{aligned} P(M \geq x) &\leq C \int_x^\infty z^{1/2r-1} e^{-A_2 z^{1/r}} \, dz \\ (5.8) \qquad &= \frac{Cr}{A_2^{1/2}} \int_{A_2 x^{1/r}}^\infty y^{-1/2} e^{-y} \, dy \\ &\leq C e^{-A_2 x^{1/r}}. \end{aligned}$$

Fix a positive integer $L$. It follows from Proposition 1.4 and (5.8) that there is a constant $C$ such that

$$\begin{aligned} \lim_{N \to \infty} \sum_{k=2}^{L} P(R_{k,N} > xN^{1-r}) &= r \sum_{k=2}^{L} P(MB_k^{1-r} > x) \\ &= r \sum_{k=2}^{L} P\Big(M > \frac{x}{B_k^{1-r}}\Big) \\ &\leq Cr \sum_{k=2}^{L} E[e^{-A_2(x/B_k^{1-r})^{1/r}}] \end{aligned}$$

$$= Cr \sum_{k=2}^{L} \int_0^1 (k-1)(1-y)^{k-2} e^{-A_2(x/y^{1-r})^{1/r}} \, dy.$$

Note that

$$\sum_{k=2}^{\infty} (k-1)(1-y)^{k-2} = \sum_{k=2}^{\infty} \sum_{i=1}^{k-1} (1-y)^{k-2}$$

$$= \sum_{i=1}^{\infty} \sum_{k=i+1}^{\infty} (1-y)^{k-2} = \sum_{i=1}^{\infty} \frac{(1-y)^{i-1}}{y} = \frac{1}{y^2}.$$

Therefore, making the substitution $z = A_2(x/y^{1-r})^{1/r}$ and using the fact that, for all real numbers $b$, there is a $C > 0$ such that $z^b e^{-z} \leq C e^{-z/2}$ for all $z \geq A_2$, we get

$$(5.9) \quad \begin{aligned} & \lim_{N \to \infty} \sum_{k=2}^{L} P(R_{k,N} > x N^{1-r}) \\ & \leq Cr \int_0^1 y^{-2} e^{-A_2(x/y^{1-r})^{1/r}} \, dy \\ & = \frac{Cr^2}{(1-r)A_2^{r/(1-r)} x^{1/(1-r)}} \int_{A_2 x^{1/r}}^{\infty} z^{(2r-1)/(1-r)} e^{-z} \, dz \\ & \leq C e^{-A_2 x^{1/r}/2}. \end{aligned}$$

By combining (5.9) with (5.8) for the $k = 1$ case, we get

$$\lim_{N \to \infty} \sum_{k=1}^{L} P(R_{k,N} > x N^{1-r}) \leq \tfrac{1}{2} C_1 e^{-C_2 x^{1/r}},$$

where $C_1$ and $C_2$ are constants that do not depend on $L$. The proposition follows by letting $\varepsilon = \tfrac{1}{2} C_1 e^{-C_2 x^{1/r}}$ and choosing $L$ as in Lemma 5.1 such that (5.2) holds for sufficiently large $N$. $\square$

**Acknowledgments.** The authors thank Anton Wakolbinger and Wolfgang Angerer for helpful discussions, and Christian Borgs and Jennifer Chayes for pointing out the connections with the preferential attachment literature. They also thank two referees for suggestions.

DEPARTMENT OF MATHEMATICS
CORNELL UNIVERSITY
ITHACA, NEW YORK 14853
USA
E-MAIL: rtd1@cornell.edu
URL: www.math.cornell.edu/˜durrett/

DEPARTMENT OF MATHEMATICS
UNIVERSITY OF CALIFORNIA, SAN DIEGO
9500 GILMAN DRIVE
LA JOLLA, CALIFORNIA 92093
USA
E-MAIL: jschwein@math.ucsd.edu
URL: www.math.ucsd.edu/˜jschwein/